\documentclass{article}

\usepackage{arxiv}

\usepackage[utf8]{inputenc} % allow utf-8 input
\usepackage[T1]{fontenc}    % use 8-bit T1 fonts
\usepackage{hyperref}       % hyperlinks
\usepackage{url}            % simple URL typesetting
\usepackage{booktabs}       % professional-quality tables
\usepackage{amsfonts}       % blackboard math symbols
\usepackage{nicefrac}       % compact symbols for 1/2, etc.
\usepackage{microtype}      % microtypography
\usepackage{lipsum}
\usepackage{graphicx}
\usepackage{caption}
\usepackage{subcaption}
\usepackage{hyperref}
\usepackage{amssymb}
\graphicspath{ {./images/} }

\title{HALE multidisciplinary ecodesign optimization with material selection}

\author{
 Edouard Duriez \\
   ICA, Université de Toulouse, ISAE-SUPAERO, MINES ALBI, UPS, INSA, CNRS\\
   3 Rue Caroline Aigle \\
  31400 Toulouse (France) \\
  \texttt{edouard.duriez@isae-supaero.fr} \\
   \And
Víctor Manuel Guadaño Martín \\
   Université de Toulouse, ISAE-SUPAERO, \\
   10 avenue edouard Belin \\
  31400 Toulouse (France) \\
  \texttt{edouard.duriez@isae-supaero.fr} \\
  \And
Joseph Morlier \\
   ICA, Université de Toulouse, ISAE-SUPAERO, MINES ALBI, UPS, INSA, CNRS\\
   3 Rue Caroline Aigle \\
  31400 Toulouse (France) \\
  \texttt{edouard.duriez@isae-supaero.fr} \\
}

\begin{document}
\maketitle
\begin{abstract}
Multidisciplinary Design Optimization (MDO) makes it possible to reach a better solution than by optimizing each discipline independently. In particular, the optimal structure of a drone won’t be the same depending on the material used. The CO2 footprint of a solar-powered High Altitude Long Endurance (HALE) drone is optimized here, the structural materials being one of the design variables. The optimization is preformed using a modified version of OpenAeroStruct, a framework based on OpenMDAO. The originality of this work is to include material choice from a discrete catalogue in the MDO approach. This is achieved through a continuous variable, enabling the use of continuous optimization algorithms. 
Our results show that, in our case, the optimal material in terms of CO2 footprint is also the optimal material in terms of weight. In order for a ``green'' surrogate material to enable a lower drone CO2 footprint, it must be almost as good in terms of weight.

\end{abstract}

% keywords can be removed

\section{Introduction}

High Altitude Long Endurance (HALE) drones driven by solar power could be an alternative to satellites for some missions. Their solar power and batteries enable them to fly for a few years, which added to their high altitude (above 20 km) make them fit for missions similar to those of satellites (\cite{dinc_performance_2021}). HALE drones can offer permanent coverage of a point or be re-positioned, and are repairable, unlike satellites, which are in orbit. Their lower altitude can offer better resolution for earth observation, but also results in smaller coverage. Their biggest advantage is their lower cost compared to satellites. HALE drones could also be more environmentally friendly by not needing a high energy consuming launcher. This advantage can be enhanced if a special attention is given to their environmental impact. For a fully electric HALE drone, most of this impact comes from the materials used and the building of the drone.

Multidisciplinary ecodesign optimization has been successfully applied to commercial planes (\cite{antoine_framework_2005}). However, HALE drone mission requirements are very different from commercial planes. Additionally to the previously mentioned features, flight speed is not a requirement as they fly in closed trajectories. Moreover a large wing surface is required in order to gather sufficient solar power. These requirements lead to very different designs, with high aspect ratios. The design of HALE drones has been studied largely already. Global configuration optimizations have been made (\cite{noth_design_2008}, \cite{morrisey_multidisciplinary_2009}, \cite{montagnier_optimisation_2010}). Scale-sized prototypes were manufactured, for mechanical testing (\cite{romeo_heliplat:_2004}, \cite{frulla_design_2008}) or aeroelastic study (\cite{cesnik_x-hale:_2012}, \cite{jones_preliminary_2015}). A multi-fidelity vehicle optimization environment (\cite{macdonald_suave_2017}) was applied to a HALE design. Optimizations seeking higher fidelity were carried out (\cite{colas_hale_2018}), and a design framework was built (\cite{hwang_solar_2019}). However, all these studies focus only on mass optimization and do not consider the environmental footprint of the drones. The aim of this work is to close this gap in a global, low-fidelity but fast, environmental impact optimization.

This work is therefore based on OpenAeroStruct (OAS), a global low-fidelity aerostructural optimization framework (\cite{jasa_open-source_2018}, \cite{chauhan_low-fidelity_2019}). This framework is itself based on OpenMDAO (\cite{gray_openmdao:_2010}), a multidisciplinary design optimization framework.

For a conventional aircraft, most of the CO2 emissions come from the fuel burned. For a solar powered HALE drone, however, no fuel is burned and most of the CO2 emissions come from the material used and the building of the drone. If the structural material used is fixed, optimizing the environmental impact of the HALE drone is equivalent to optimizing the mass of the drone, as detailed in section \ref{mat}. Therefore, we decide to include the choice of the structural material in the optimization. A recurring problem in optimal material selection is that the material variable is a discrete variable. Another problem is that optimal properties are usually antagonist. For example a material with a higher Young's modulus will also have a higher density.

In order to counter these difficulties Ashby indexes can be used (\cite{ashby_materials_2004}). These enable to order materials according to a unique index. However, this method only works for choosing a material for a given simple part under a given loading. This method has been extended for simultaneous material selection and geometry design (\cite{rakshit_simultaneous_2007}). Yet, this method is limited to a structure and cannot take into account the impact of this structure on a bigger system and the resulting feedback loops. Making the material variable continuous has been proposed to find the optimal material (\cite{bendsoe_analytical_1994}) and developed in the form of an interpolation in topology optimization (\cite{zuo_multi-material_2017}). However, these methods are not used in multidisciplinary design optimizations. One of the main contributions of this work is to add the discrete choice of material in a continuous multidisciplinary optimization framework.

Having already given a description of related works, the rest of this paper is organized as follows. Section \ref{sec_mdo} presents the HALE drone model and optimization framework where the material choice is a variable. Section \ref{secresults} validates this framework and details the obtained results and how to interpret them. Finally, Section \ref{conclu} gives some concluding remarks.

\section{HALE drone MDO formulation}
\label{sec_mdo}
\subsection{OpenAeroStruct to Eco-HALE}
\label{OAStoEH}
\begin{sloppypar}
We choose to derive our work from the OpenAeroStruct with wingbox framework (\cite{chauhan_low-fidelity_2019}). OpenAeroStruct is originally made to study commercial aircraft. In order to model a HALE drone, some changes have to be made.
As we choose to use a comparison with the single-boom HALE drone from \cite{colas_hale_2018} to validate our framework, we use the same data as in this paper each time it is possible.
\end{sloppypar}

First of all, as there is no need for a 2.5G maneuver on a HALE drone, we change the design points used. We use one cruise flight design point where the lift and power constraints must be satisfied, and one gust design point where the structural constraints must be satisfied.

\begin{sloppypar}
The OpenAeroStruct constraint on performance, based on the Breguet equation cannot be used for a solar-powered HALE aircraft, as the total weight of the aircraft does not change during flight and the range is not limited. Therefore we suppress it, and also suppress all components regarding fuel.
Instead, it is a power equilibrium that gives a performance constraint.
Indeed, the power used by the propulsion and the payload must be produced by the solar panels or stored in the batteries.
The power needed for propulsion $P_{prop}$ is derived from the 1D equilibrium between thrust ($T$) and drag ($D$).
Thrust and drag are expressed in Eqs. \ref{eqDrag} and \ref{eqThrust}, where $W$ is the total weight of the drone, $C_d$ is its drag coefficient, $C_l$ is its lift coefficient, and $v$ its speed. 
\end{sloppypar}
\begin{equation}
\label{eqDrag}
D = W \cdot \frac{C_d}{C_l}
\end{equation}
\begin{equation}
\label{eqThrust}
P_{prop} = T \cdot v
\end{equation}

These give a link between the power needed for propulsion and the total weight of the drone: Eq. \ref{eqProp}
\begin{equation}
\label{eqProp}
P_{prop} = W \cdot v \cdot \frac{C_d}{C_l}
\end{equation}

The relationship we actually use is Eq. \ref{eqProp2}, in order to account for the propulsion efficiency $\eta$ and the power needed by the payload and the avionics $P_{payload}$\label{neededP}.

\begin{equation}
\label{eqProp2}
P_{needed} = \frac{W \cdot v \cdot C_d}{C_l \cdot \eta} + P_{payload}
\end{equation}

During day time, the power is supplied by the solar cells. Therefore a minimum wing surface ($S_{wing}$) is needed in order to place a sufficient number of solar cells. This is expressed in Ineq. \ref{eqSurf}.
\begin{equation}
\label{eqSurf}
S_{wing} > \frac{P_{needed}}{A_{PV}}
\end{equation}

\begin{figure}
  \centering
  \includegraphics[width=0.49\textwidth]{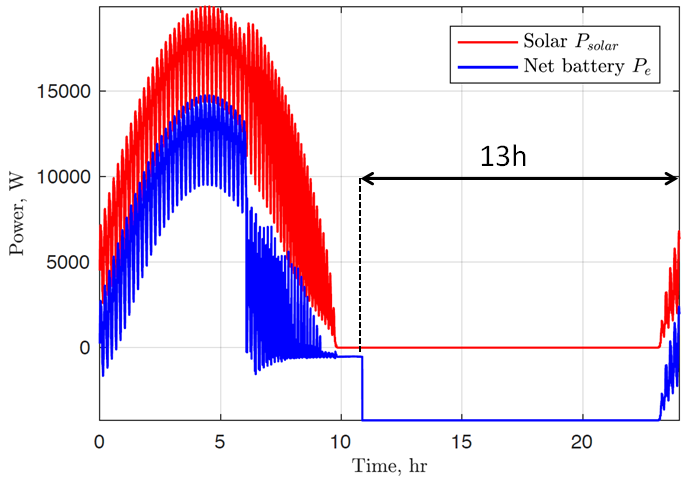}
\caption{Solar and battery power for FBHALE (\cite{colas_hale_2018}). The night time we consider starts at the end of the gliding phase.}
\label{powerFBH}       % Give a unique label
\end{figure}

\begin{sloppypar}
$A_{PV}$ is the power produced per unit of area of solar cells. Ineq. \ref{eqSurf} is the constraint on performance we use.
$A_{PV}$  is chosen carefully in order to be in the same case as FBHALE (\cite{colas_hale_2018}). Fig. \ref{powerFBH} comes from \cite{colas_hale_2018}, it was used to extract the mean solar power harvested during the 24h power equilibrium considered. This mean power was divided by the area of solar cells in \cite{colas_hale_2018}, and leads to $A_{PV}=54W/m^2$. 
\end{sloppypar}

\begin{figure}
  \centering
  \includegraphics[width=0.49\textwidth]{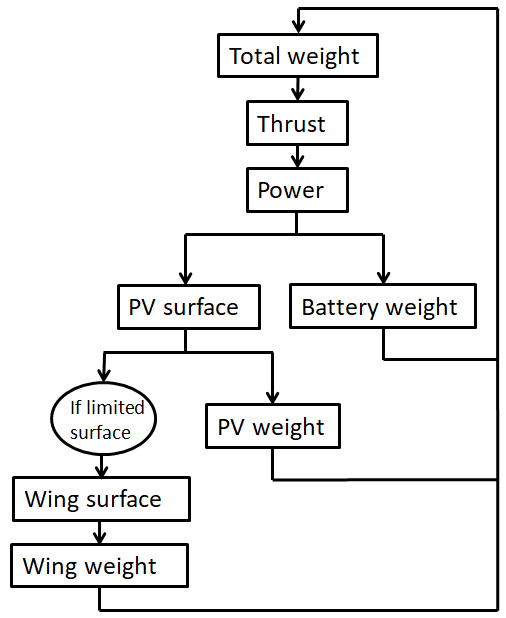}
\caption{Power and weight coupling: if the total weight changes, the power needed to thrust that weight changes; if the needed power changes, the battery weight, the photo-voltaic cells (PV) weight and possibly the wing weight, change, resulting in a total weight change. This coupling is solved by the same Gauss-Seidel solver as the aerostructural coupling.}
\label{powerWeight}       % Give a unique label
\end{figure}

\begin{figure}
  \centering
  \includegraphics[width=0.49\textwidth]{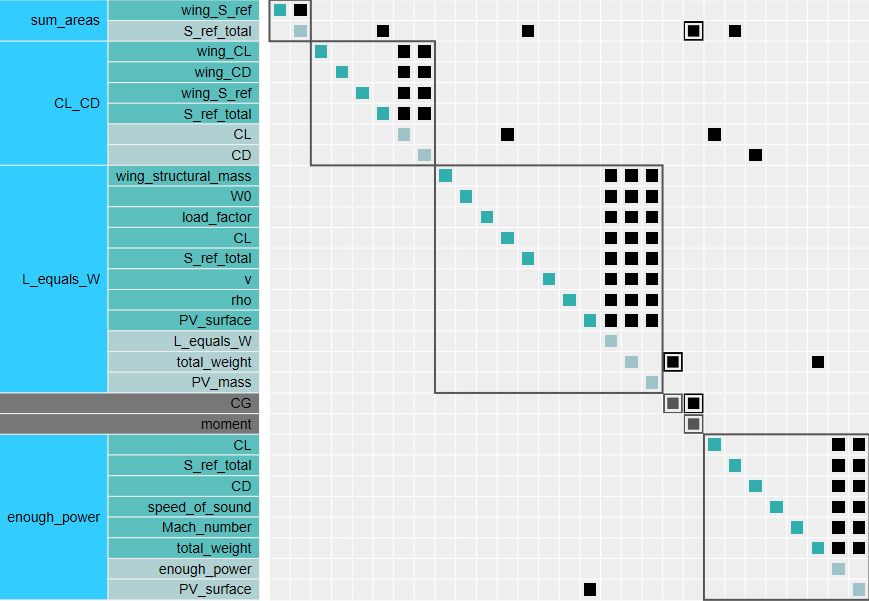}
\caption{Component layout of the performance group adapted from OpenAeroStruct \cite{jasa_open-source_2018}, \cite{chauhan_low-fidelity_2019}. The coupling between power and weight is visible as a “PV\_surface” variable from which the weight of the batteries and solar panels are derived.}
\label{schema}       % Give a unique label
\end{figure}

The computed surface of solar cells necessary to produce the needed power is taken into account in the total weight computation, meaning that weight and power are coupled, as can be seen on Figs. \ref{powerWeight} and \ref{schema}.

\begin{sloppypar}
The mass of the solar cells added to the total mass is computed using a surface density from \cite{colas_hale_2018}. 
\end{sloppypar}

The mass of the batteries is also taken into account. In order to compute this mass, an assumption that the batteries power the drone during a 13h night time ($t_{night}$) is made. This assumption is also based on Fig. \ref{powerFBH}, and enables a great simplification of the problem by not modeling the flight and gliding phase. This time multiplied by the power needed gives the size of the battery in kWh\label{batnight}, which in turn gives the mass of the batteries ($M_{bat}$) when divided by the energy density ($d_{bat}$ in kWh/kg). The energy density was taken equal to that of FBHALE. The mass of the batteries is therefore computed with Eq. \ref{eqBat}. Taking the battery mass into account adds to the weight and power coupling, as on Fig. \ref{powerWeight}.

\begin{equation}
\label{eqBat}
M_{bat}=\frac{P_{needed} \cdot t_{night}}{d_{bat}}
\end{equation}

These masses of the batteries and the solar panels are distributed along each of the beam elements of the wing, as was done in OpenAeroStruct for the fuel mass. This is closer to reality than a concentrated mass, and relieves the wing enabling thinner structural skins. Other masses are also taken into account but not distributed, such as the avionics and payload mass, considered as fixed during optimization, and the propulsion system mass.

The propulsion mass ($M_{prop}$) models both the mass of the engine and the propeller. It is obtained as the product of $P_{prop}$ (Eq. \ref{eqProp}) and the propulsion density ($d_{prop}$ in kg/W), as in Eq. \ref{mprop}. In order for the data to be as close as possible to FBHALE to validate our framework, without modelling the propeller, this propulsion density is estimated based on \cite{colas_hale_2018}. This estimate is obtained by dividing the propulsion mass (extracted from Fig. \ref{fbmasscomp}) by the power used for propulsion (extracted from Fig. \ref{powerFBH}). These engines and propellers are added as point masses onto the wing structure. As in Eq. \ref{mmot}, the mass of each motor ($M_{mot}$) is obtained by dividing the propulsion mass by the number of motors ($n_{mot}$). \label{propulsion}

\begin{equation}
\label{mprop}
M_{prop}=P_{prop} \cdot d_{prop}
\end{equation}

\begin{equation}
\label{mmot}
M_{mot}=\frac{M_{prop}}{n_{mot}}
\end{equation}

The power system (solar panels, maximum power point tracker and batteries) is sized for a worst case corresponding to the winter solstice where the least solar power is available. By choosing a more favourable launch day, much more power is available during climb. Therefore, the climb phase does not influence the sizing of the power system. As a result, computing the propulsion mass as in Eq. \ref{mprop} also enables us not to model the climb phase of the drone. This phase is indeed taken into account through the propulsion mass data from \cite{colas_hale_2018}. 

A shear gust wall is added to the framework, in order to compute the loads used for the sizing of the wings. It has the same magnitude as in \cite{colas_hale_2018}.

Our model only takes into account the weight of the wing, the batteries, the solar panels, the propulsion, the maximum power point tracker (MPPT), the avionics (only the fixed mass) and the payload mass. In order to account for the extra mass of the harness, the landing gear, the interfaces, the horizontal and vertical tail, the boom and the pod, the total weight of the drone is increased by 10\%. This value comes from Fig. \ref{fbmasscomp} (\cite{colas_hale_2018}).

\begin{figure}
  \centering
  \includegraphics[width=0.49\textwidth]{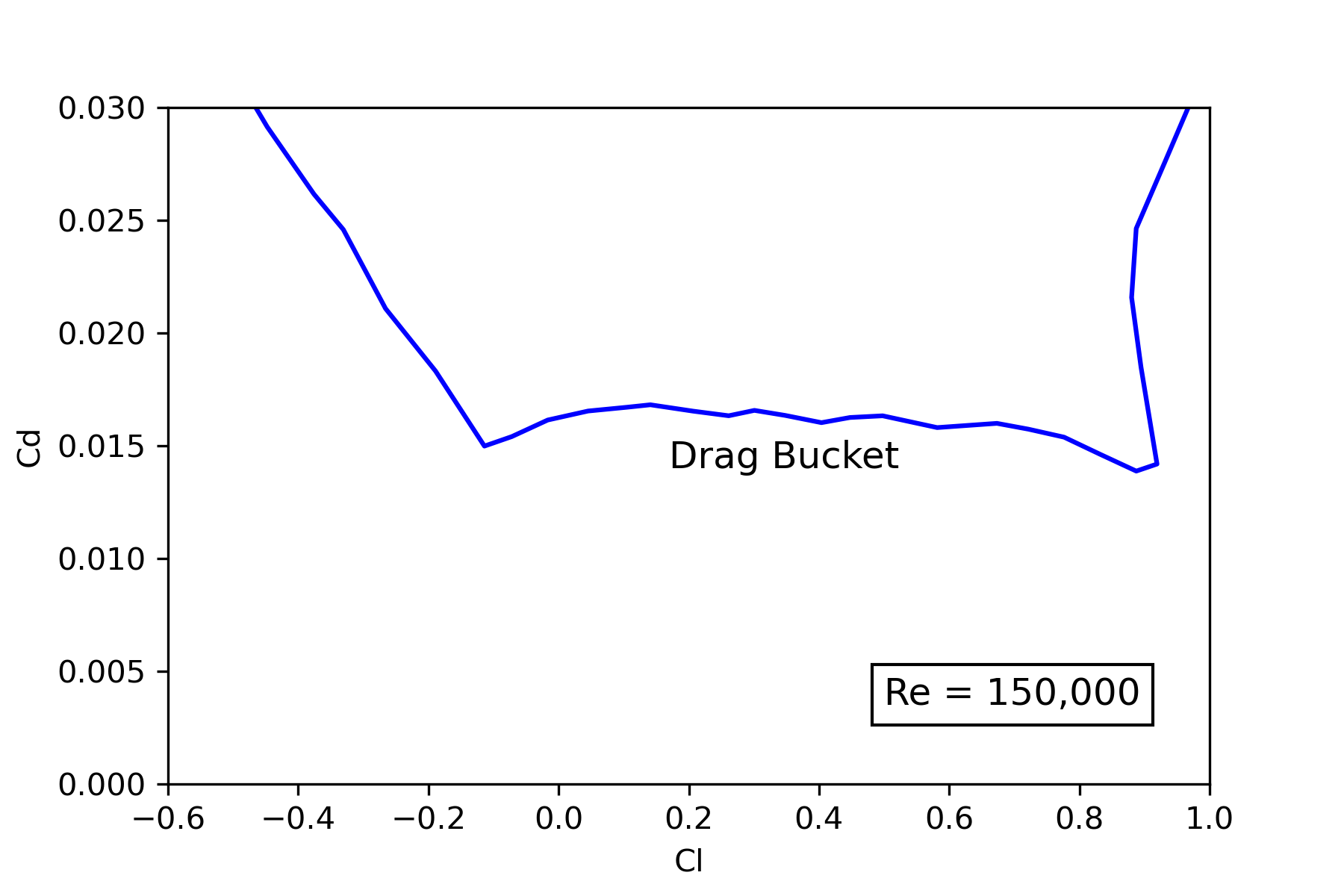}
\caption{Polar curve for NACA 63412 airfoil.}
\label{polar}       % Give a unique label
\end{figure}

Finally, the operational conditions of the HALE, such as the high altitude and the low speed, imply a low Reynolds number that lies in the 150,000 – 200,000 range. For this reason, high values of lift coefficient ($C_l$) are needed and a low-drag configuration is mandatory in order to minimize the power needed for propulsion (Eq. \ref{eqProp}). Thus, as proposed in \cite{mattos2013optimal}, the NACA 63412 profile is selected. It is a laminar airfoil with high maximum lift coefficient ($Cl_{max}$), low value of the moment coefficient ($C_m$) and ``drag bucket'' covering the desired $C_l$ range (see Fig. \ref{polar}).

\subsection{MDA framework summary}

The multidisciplinary design analysis (MDA) framework is summarized in Tab. \ref{disc}. The method modeling each discipline is shown in this table, as well as the implementation used and a reference. We use Gauss-Seidel fixed-point iterations to converge the MDA. This implies that the last outputs from the other analysis are used to run each analysis, until MDA convergence.

\begin{table*}
\centering
\caption{Summary of multidisciplinary design analysis framework.}
\label{disc}       % Give a unique label
\begin{tabular}{llll}
\hline\noalign{\smallskip}
Discipline & Method & Implementation & Reference\\
\noalign{\smallskip}\hline\noalign{\smallskip}
  Aerodynamics & VLM & OAS & \cite{anderson_fundamentals_1991} \\
  Structure & Wingbox beams & OAS & \cite{chauhan_low-fidelity_2019} \\
  Energy & Simple in-house method & Section \ref{OAStoEH} & data from \cite{colas_hale_2018} \\
  Environmental & Proportional to mass & Section \ref{objfunc} & data from \cite{wetzel_update_2015}, \cite{hao_ghg_2017} \\
\noalign{\smallskip}\hline
\end{tabular}
\end{table*}

The energy and environmental models are described in other sections of this paper, listed in Tab. \ref{disc}. We briefly present hereafter the aerodynamics and structural models used. They are kept unchanged from the OpenAeroStruct (\cite{jasa_open-source_2018}) implementation, apart from the changes to the structure model described in Section \ref{addCons}.

The aerodynamics model is the implementation of the VLM model (\cite{anderson_fundamentals_1991}) by OpenAeroStruct (\cite{jasa_open-source_2018}). Each panel of the lifting surface is modeled as a horseshoe vortex. This horseshoe vortex is made of three vortex segments : one span-wise bound vortex located at $\frac{1}{4}^{th}$ of the panel's chord-wise length and two semi-infinite chord-wise vortices extending towards infinity in the free-stream direction. The total lifting surface is modeled as the superimposition of all these horseshoe vortices. Flow tangency conditions are imposed preventing normal flow through the panel. The aerodynamic forces on each panel are derived from the circulation strength of the vortices. Finally, lift and drag are derived from these forces. Skin friction drag is also estimated. The reader is referred to \cite{jasa_open-source_2018} for more detail on this implementation.

\begin{sloppypar}
The structure model is the beam FEM with wing-box cross-sections implemented in OpenAeroStruct (\cite{chauhan_low-fidelity_2019}). The span-wise discretization of this model (7 elements) is the same as for the VLM mesh. Each span-wise section of wing is modeled as a beam element. The cross section of this beam element is a hollow rectangle modeling the upper and lower skins of the wing and the two spars. This model is used to compute the worse case von Mises stress values that is used as a mechanical failure constraint. The reader is referred to \cite{chauhan_low-fidelity_2019} for more detail on this implementation.
\end{sloppypar}

\subsection{MDO framework}

\subsubsection{Additional constraints}
\label{addCons}

Having only a mechanical failure constraint leads to too thin skin thicknesses. Therefore, a buckling constraint is added.

This buckling constraint imposes that the stress in the top skin of each element stays below the buckling critical stress of that panel. In order to have only one constraint, a Kreisselmeier-Steinhauser (KS) aggregation is used (\cite{kreisselmeier_systematic_1979}, \cite{lambe_evaluation_2017}), as for the initial mechanical failure constraint.

The buckling critical stresses are estimated by the case of a homogeneous rectangular curved plate subject to combined axial compression and shear. Therefore, two different buckling critical stresses are considered: an axial critical stress ($\sigma_c$) and a shear critical stress ($\tau_c$), which can be put in the form of Eq. \ref{bucaxial} and Eq. \ref{bucshear}, respectively, where $k_c$ is the buckling coefficient for axial load, $k_s$ the buckling coefficient for shear, $b$ the width of the panel, and $D$, the flexural stiffness. In our approximate buckling case, we take $b$ to be the distance between the two spars (fixed equal to half the chord) and the values of $k_c$ and $k_s$ can be obtained graphically from the charts presented in \cite{gerard1957handbook}, using the geometry of the panel as an input. \label{bucsection}

\begin{equation}
\label{bucaxial}
\sigma_c=k_c \cdot \frac{\pi^2 \cdot D}{b^2}
\end{equation}

\begin{equation}
\label{bucshear}
\tau_c=k_s \cdot \frac{\pi^2 \cdot D}{b^2}
\end{equation}

Thus, the buckling constraint is approximated by a parabolic interaction, as in Ineq. \ref{bucinter}, proposed in \cite{gerard1957handbook}, $R_s$ and $R_c$ being the stress ratios for shear and axial compression, respectively. These stress ratios are defined as the ratio of the stress at buckling under combined loading to
the buckling stress under simple loading.

\begin{equation}
\label{bucinter}
R_s^2 + R_c < 1
\end{equation}

We also add a constraint imposing that the skin thickness cannot be greater than half the wing thickness, in order to prevent the top and bottom skin from intersecting at the wing tips, where the chord is smaller and the thickness-to-chord ratio of the wing is small too.

\subsubsection{Objective function}
\label{objfunc}

We choose the CO2 emitted by the HALE drone during its life cycle ($CO2_{tot}$) as our objective function. As the drone does not burn any fuel during flight, it is assumed that the CO2 is mainly emitted before the use of the drone and is in particular due to the materials used and their processing. As the aim is to have a simple and fast tool, we consider only the CO2 emitted by the material used for the structure ($CO2_{struct}$), the CO2 emitted by the solar panels ($CO2_{PV}$), and the CO2 emitted by the batteries ($CO2_{bat}$), as these have the most influence on the total emitted CO2. This is shown on Eq. \ref{eqCOtot}.

\begin{equation}
\label{eqCOtot}
CO2_{tot} = CO2_{struct} + CO2_{PV} + CO2_{bat}
\end{equation}

The CO2 emitted by the structure is computed as the product of the mass of the spars ($M_{spar}$) by the CO2 footprint of the material used for the spars ($CO2_{mat1}$) and the product of the mass of the skins ($M_{skin}$) by the CO2 footprint of the material used for the skins ($CO2_{mat2}$), as in Eq. \ref{eqCOstr}. This CO2 footprint of the materials is considered to be equal to a weighted sum of the CO2 footprint of their primary production ($CO2_p$) and the CO2 footprint of their recycling ($CO2_r$), depending on their recycled fraction in current supplies ($\eta_r$), as in Eq. \ref{eqCOmat}.

\begin{equation}
\label{eqCOstr}
CO2_{struct} = M_{spar} \cdot CO2_{mat1} + M_{skin} \cdot CO2_{mat2}
\end{equation}
\begin{equation}
\label{eqCOmat}
CO2_{mat} = \eta_{r} \cdot CO2_{r} + (1-\eta_{r}) \cdot CO2_{p}
\end{equation}

It can be seen in Eqs. \ref{eqCOstr} and \ref{eqCOmat}, that CO2 footprint of the processing of the materials is not taken into account. We make this choice because CO2 footprint of processing varies a lot depending on the process, and we have no information related to this stage. This choice is acceptable as the footprint of the most probable processes are much lower than the footprint of the material primary production, and this is especially true for the materials that were found to be optimal.

The CO2 emitted by the solar panels is computed as the product of the power needed on board the HALE ($P_{needed}$ detailed in section \ref{neededP}) and the emissions per power ($CO2_{/W}$), as in Eq. \ref{eqCOpv}. The value of $CO2_{/W}$ is taken from \cite{wetzel_update_2015}.

\begin{equation}
\label{eqCOpv}
CO2_{PV} = P_{needed} \cdot CO2_{/W}
\end{equation}

The CO2 emitted by the batteries is computed as the product of the size of the battery ($P_{needed} \cdot t_{night}$ as in section \ref{batnight}) and the emissions per Wh ($CO2_{/Wh}$), as in Eq. \ref{eqCObat}. The batteries are assumed to be close to electric car Li-ion batteries, and the value of the emissions per Wh is taken equal to the mean of the three batteries studied in \cite{hao_ghg_2017}.

\begin{equation}
\label{eqCObat}
CO2_{bat}=P_{needed} \cdot t_{night} \cdot CO2_{/Wh}
\end{equation}

\subsubsection{Design variables}

We use eight geometric design variables and a material design variable which we discuss in detail in section \ref{mat}. B-splines parameterized using four control points each are adopted in order to vary some geometric design variables along the semi-span.

The only angle design variable we use is the vector of control points for the twist, because it also enables to control the angle of attack. Therefore, our variable twist represents the sum of geometric twist and angle of attack at cruise conditions.

\begin{sloppypar}
We keep the wingbox design variables of OpenAeroStruct, which are the control points for the skin thickness and the spar thickness. They are mainly determined by the buckling constraint and the mechanical failure constraint, respectively. We use the same thickness distribution for both the upper and lower skins, and the same thickness distribution for both the forward and rear spars.
\end{sloppypar}

For the wing geometry, we keep the control points for the thickness-to-chord ratio of the wing as a design variable, and add the span, chord and taper ratio as design variables, in order to be able to satisfy the power constraint and the weight-equals-lift constraint. These variables are able to generate the particular geometries of HALE drones with very high aspect ratios.

Additionally, the motor spanwise location is added as a new design variable. Since we are considering a symetrical twin-motor HALE, this new design variable is defined as the ratio of the distance between the plane of symmetry and the motor to the semi-span of the wing.

Some limitations have to be imposed to the design space in order to compensate for some missing physics in the model. Wings cannot be too tapered because tip stall can happen otherwise. As stall is not modelled in our framework, we choose a lower limit of 0.3 for the taper ratio design variable. For the same reason, a lower and upper twist limits of -15 and 15 degrees are considered. In this way, the critical or stalling angle of attack is taken into account. Furthermore, wings cannot be too narrow because 1-cosine gusts limit their aspect ratio. As these gusts are not modelled, we also add a lower limit of 1.4 metres for the root chord.

\subsubsection{Material formulation}
\label{mat}

The battery and solar panel types are fixed in our case. Therefore minimizing the CO2 emissions due to the solar panels and batteries means minimizing the power needed, which means minimizing the mass of the drone to be thrusted. Similarly, minimizing the CO2 emissions due to the structural materials means minimizing the mass of materials used, if the choice of materials is not a variable. Therefore, if the structural materials were fixed, optimizing the CO2 impact of the drone would be equivalent to optimizing the mass of the drone. Therefore, we choose to add the choice of the structural materials as a design variable. This additional design variable takes the form of a vector with two components: the density of the material used for the spars and the density of the material used for the skins. 

In order to have a unique material design variable, we choose to access materials through their density. The material data we need (Young's modulus, shear modulus, yield strength, and CO2 emissions) is therefore accessed as a function of the density design variable.

\begin{figure}
  \centering
  \includegraphics[width=0.49\textwidth]{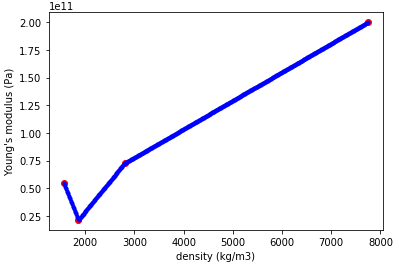}
\caption{Young’s modulus example of linear interpolation of materials.}
\label{linint}       % Give a unique label
\end{figure}

In order to have only continuous variables, we decide to make the density design variable continuous. This is achieved by interpolating linearly each material property in the space between real materials. The resulting interpolation is shown on Fig. \ref{linint} for four real materials: carbon fibre reinforced composite (CFRP), glass fibre reinforced composite (GFRP), aluminium and steel.

However, with this interpolation, it is possible to find an optimum in an interpolated material corresponding to no real material. In order to get rid of this possibility, the interpolation is penalized by the use of a power term in the interpolation, as in \cite{zuo_multi-material_2017}. An example of this interpolation for the Young’s modulus is given in Eq. \ref{interp}, where $E$ is the interpolated Young’s modulus at density $\rho$, $\rho_i$ and $\rho_{i+1}$ are the densities of the real materials framing the one being interpolated, $E_i$ and $E_{i+1}$ are the respective Young’s modulus of these real materials, and $p$ is the power used for penalization.

\begin{equation}
\label{interp}
E(\rho)=A \cdot \rho^p+B \\
\end{equation}
with $A=\frac{E_{i+1}-E_i}{\rho_{i+1}^p-\rho_i^p}$  and $B=E_i-A \cdot \rho_i^p$ \\

\begin{figure}
  \centering
  \includegraphics[width=0.49\textwidth]{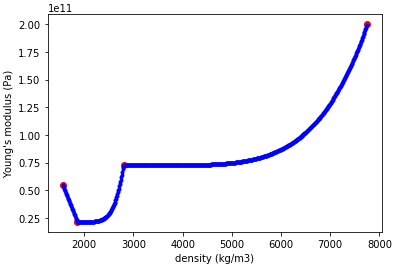}
\caption{Young's modulus example of penalized interpolation of materials.}
\label{penint}       % Give a unique label
\end{figure}

The resulting interpolation is shown on Fig. \ref{penint} for $p=5$, and for the same real materials as on Fig. \ref{linint}.

To the left of a real material, a small decrease in density brings for example a great decrease in Young’s modulus. Therefore, the slight advantage due to a slightly lower density is compensated by the Young's modulus evolution. To the right of a real material, even a big increase in density brings a very small increase in Young’s modulus. Therefore, the slight advantage due to a slightly higher Young's modulus is compensated by the higher density. As a result, choosing a material that is not real is not beneficial, and the final optimum will be a real material.

\begin{figure}
  \centering
  \includegraphics[width=0.49\textwidth]{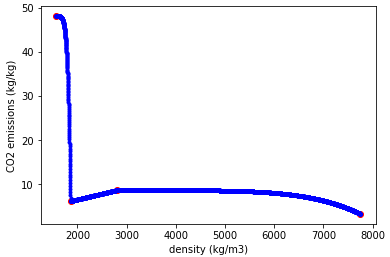}
\caption{CO2 emissions example of penalized interpolation.}
\label{negint}       % Give a unique label
\end{figure}

In order not to advantage fake materials, the penalization isn’t used in the case where a lower density brings a better property. This appears for example in the left-hand part of Fig. \ref{penint}, between CFRP and GFRP, where decrease in density brings an increase in Young’s modulus.
For material properties for which a smaller value is better, such as CO2 emissions, the penalization is inverted, as on Fig. \ref{negint}.

This penalization must be given to the interpolation only at the end of the optimization. Indeed, as we use a gradient optimization, once the penalization is given to the interpolation, the solution is stuck in the local minimum corresponding to the closest material.

The optimal material being an angular point, it is often the optimum found even without any penalization in the interpolation. However, adding the penalization is more robust.

\begin{table*}
\centering
\caption{Summary of multidisciplinary design optimization framework.}
\label{MDOsum}       % Give a unique label
\begin{tabular}{lll}
\hline\noalign{\smallskip}
\textbf{Objective function} & \textbf{Dimension}  & \textbf{Bounds} \\
\noalign{\smallskip}\hline\noalign{\smallskip}
$CO2_{tot}$ & $\mathbb{R}$  &  \\
\noalign{\smallskip}\hline\noalign{\smallskip}
\textbf{Design variables} &   &  \\
\noalign{\smallskip}\hline\noalign{\smallskip}
  Density & $\mathbb{R}^2$  & [400, 8000] kg/m$^3$ \\
  Twist control points & $\mathbb{R}^4$  & [-15, 15] deg \\
  Skin thickness ($t_{skin}$) control points & $\mathbb{R}^4$  & [0.001, 0.1] m \\
  Spar thickness control points & $\mathbb{R}^4$  & [0.001, 0.1] m \\
  Thickness-to-chord ratio control points & $\mathbb{R}^4$  & [0.01, 0.4] \\
  Span & $\mathbb{R}$  & [1, 1000] m \\
  Root chord & $\mathbb{R}$  & [1.4, 500] m \\
  Taper ratio & $\mathbb{R}$  & [0.3, 0.99] \\
  Motor location over semi-span ratio  & $\mathbb{R}$  & [0, 1] \\
\noalign{\smallskip}\hline\noalign{\smallskip}
\textbf{Constraints} &   &  \\
\noalign{\smallskip}\hline\noalign{\smallskip}
  Mechanical failure $\sigma < \sigma_{max}$ & $\mathbb{R}^7$  &  \\
  Buckling $R_s^2 + R_c < 1$ & $\mathbb{R}^7$  &  \\
  Skin thickness $2 t_{skin} < t_{wing}$ & $\mathbb{R}^4$  & \\
  Power equilibrium $P_{needed} / A_{PV} < S_{wing}$ & $\mathbb{R}$ & \\
\noalign{\smallskip}\hline
\end{tabular}
\end{table*}

\subsubsection{MDO framework summary}

Tab. \ref{MDOsum} summarizes the MDO framework. 

It is important to highlight that SLSQP (Sequential Least Squares Programming) optimization algorithm is used in this work. The chosen optimizer options (stopping criteria) are $10^{-3}$ for the convergence accuracy and 250 for the maximum number of iterations. We check the analytical derivatives using the complex step method (\cite{martins_complex-step_2003}).

\begin{figure*}
\begin{minipage}{.48\textwidth}
\begin{subfigure}{\linewidth}
  \centering
  \includegraphics[width=\linewidth]{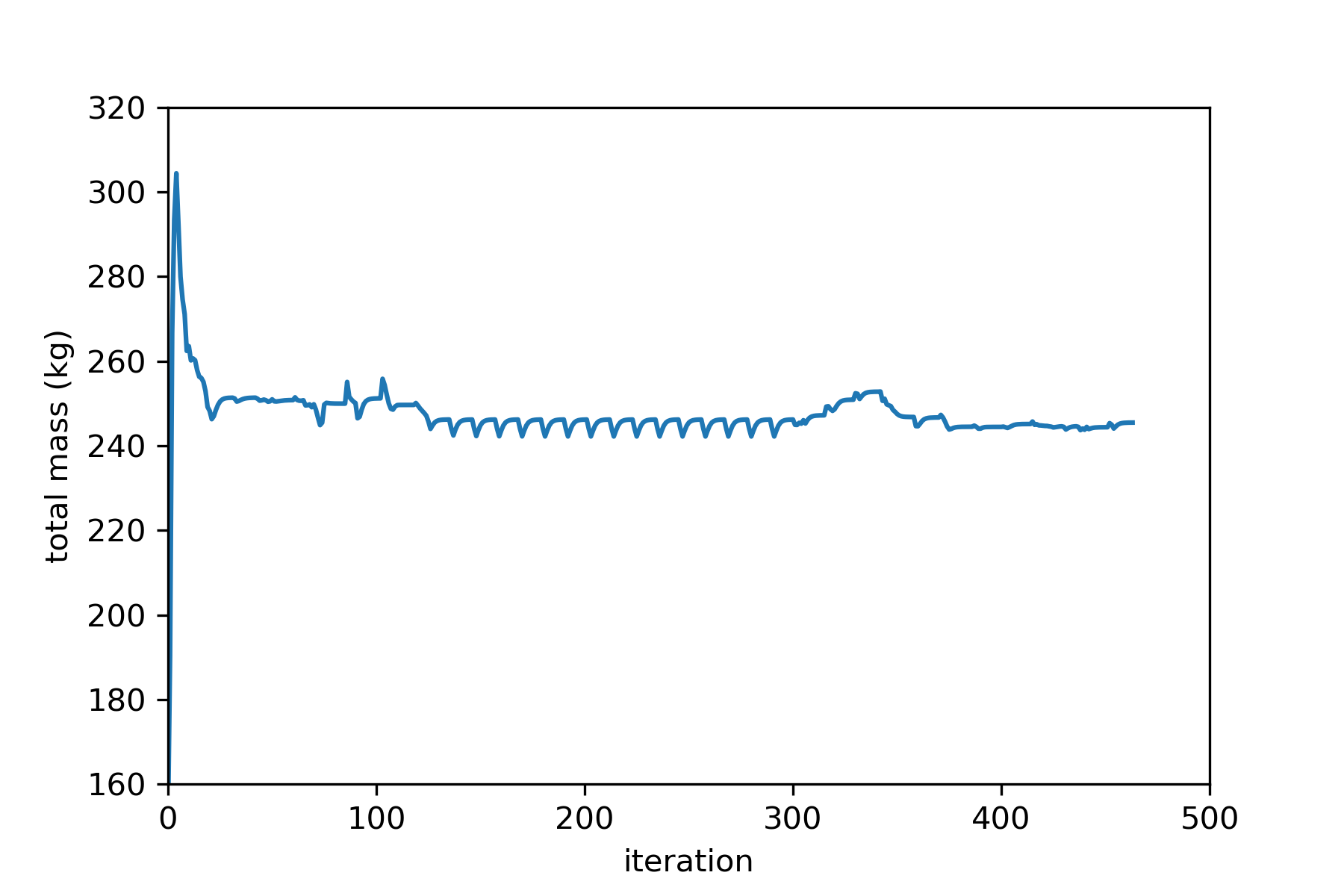}
  \caption{Objective function: total mass of the drone.}
  \label{fig:sub3}
\end{subfigure}\\[\baselineskip]
\begin{subfigure}{\linewidth}
  \vspace{-2.4mm}
  \includegraphics[width=\linewidth]{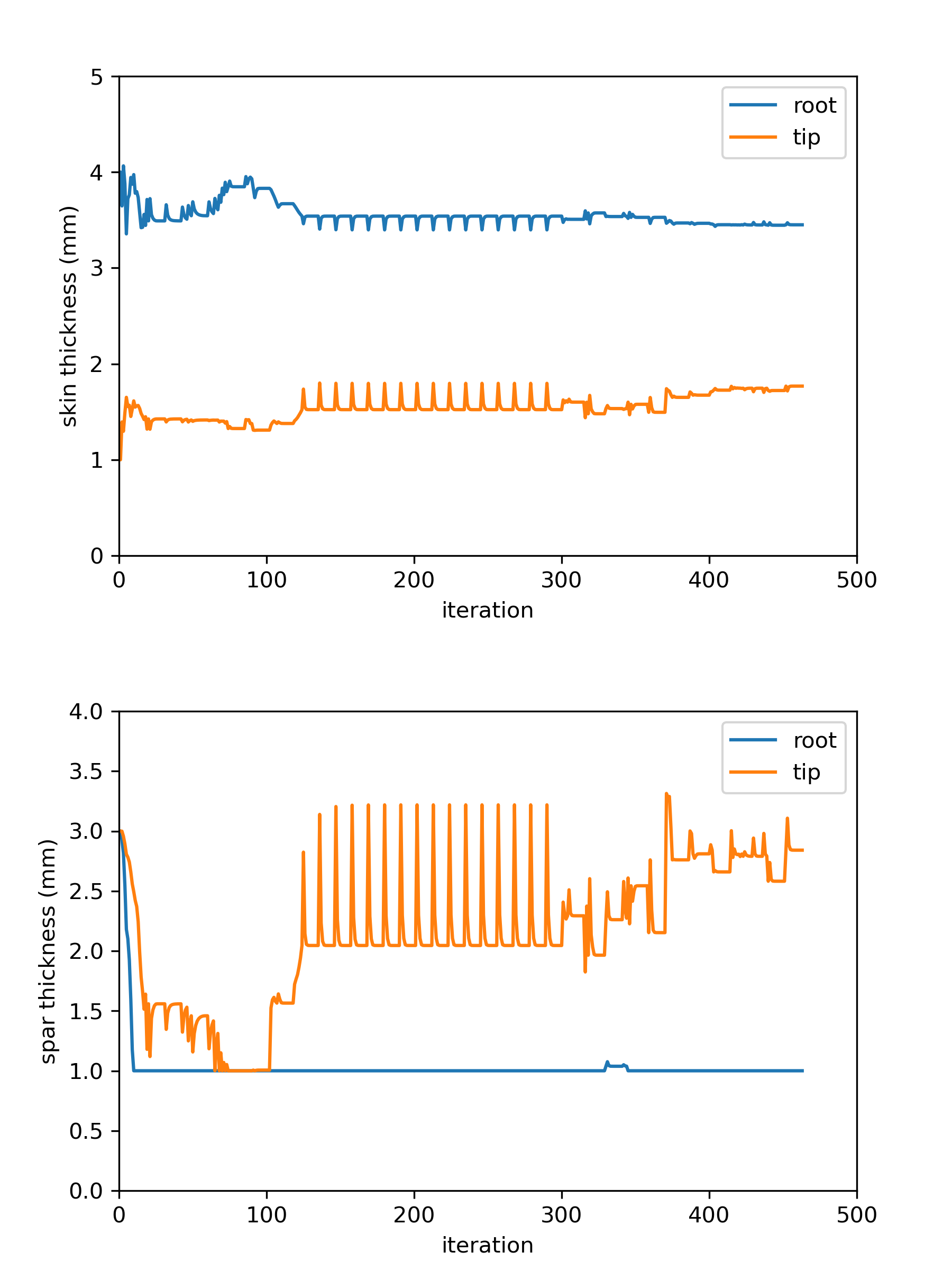}
  \caption{Skin and spar thicknesses.}
  \label{fig:sub5}
\end{subfigure}%
\end{minipage}
\begin{subfigure}{.5\textwidth}
  \centering
  \includegraphics[width=\linewidth]{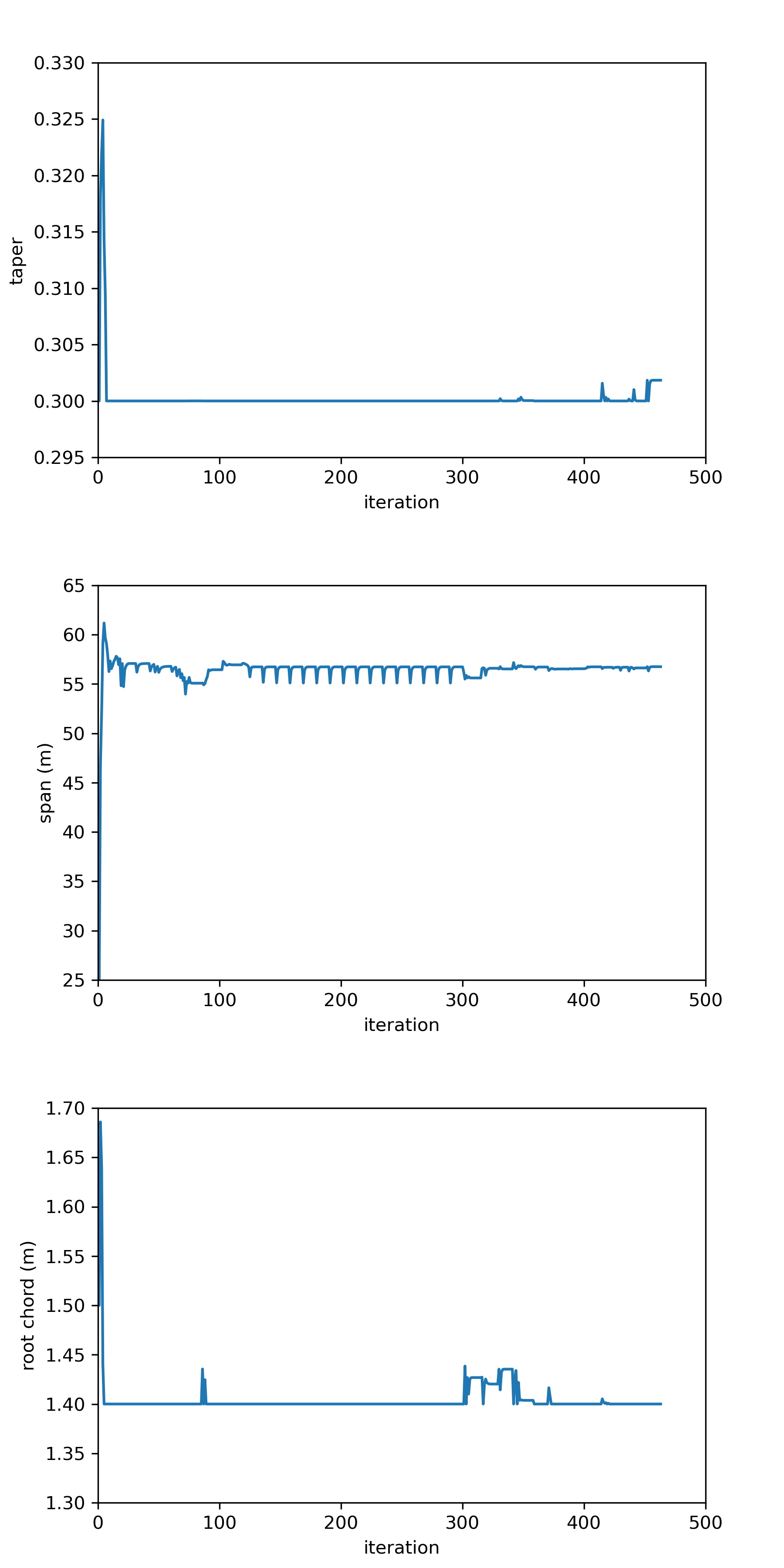}
  \caption{Geometrical variables.}
  \label{fig:sub4}
\end{subfigure}
\caption{Convergence graphs of some design variables for validation case. All final values are summarized on Tab. \ref{valid}. It can be seen that most of the design variables tend to a stable final value except for the spar thickness at the wing tip which fluctuates.}
\label{convgraph}
\end{figure*}

\section{Results}
\label{secresults}
\subsection{Mass minimization problem (validation with FBHALE)}

In order to validate our framework, we compare it to \cite{colas_hale_2018} on a similar run. In order for the comparison to be significant, we use the same data as \cite{colas_hale_2018} every time it is possible. For this validation, we also change our objective function in order to do a total mass optimization as in \cite{colas_hale_2018}, instead of a CO2 optimization. Finally, we fix the material design variable to a material similar to that of \cite{colas_hale_2018}.

\begin{figure*}
\begin{subfigure}{.49\textwidth}
  \centering
  \includegraphics[trim={0.cm 1.32cm 0.cm 1.35cm}, clip, height=5cm]{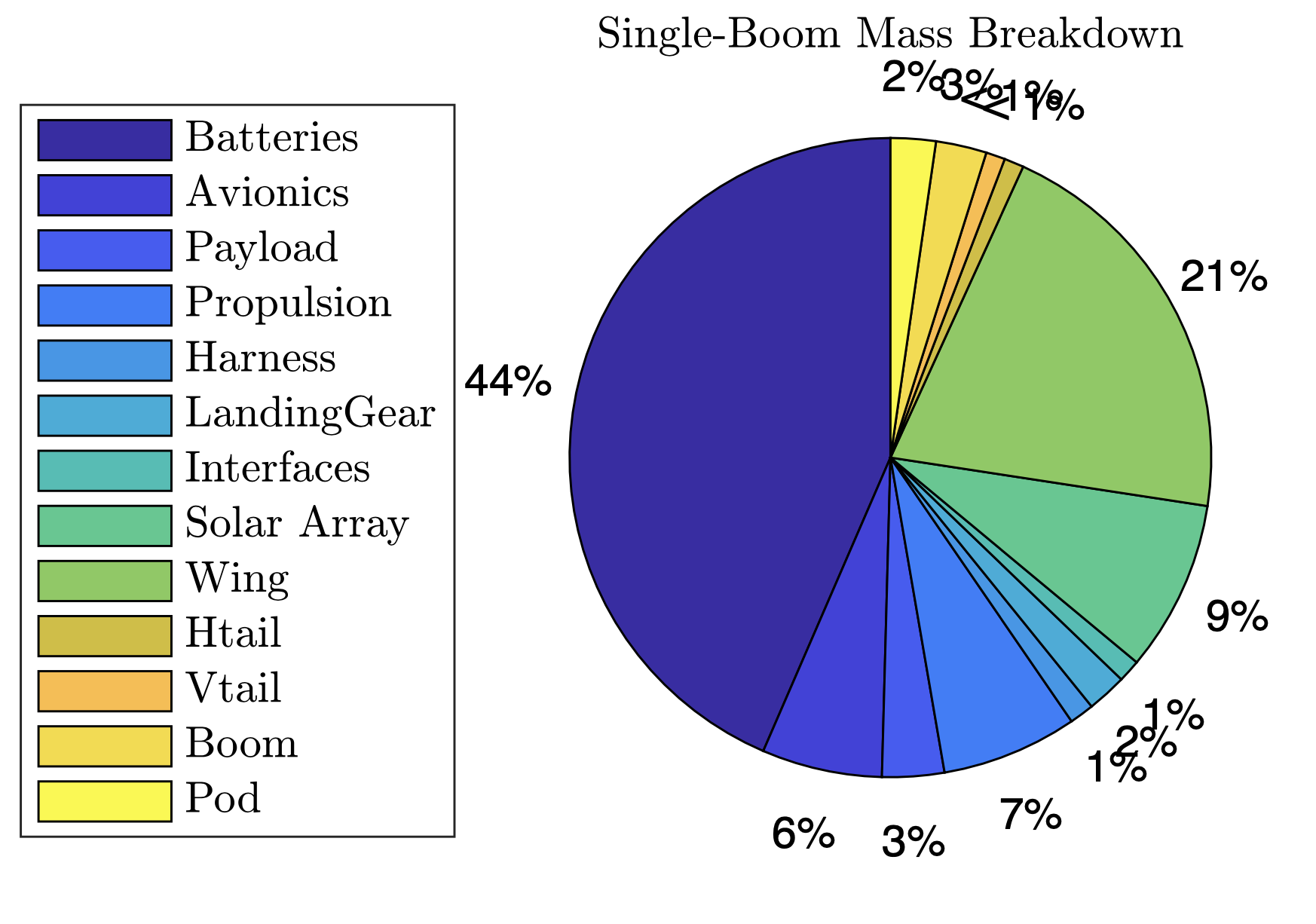}
  \caption{FBHALE’s single-boom (\cite{colas_hale_2018}).}
  \label{fbmasscomp}
\end{subfigure}
\begin{subfigure}{.49\textwidth}
  \centering
  \includegraphics[trim={8cm 1.5cm 1.5cm 1.5cm}, clip, height=5cm]{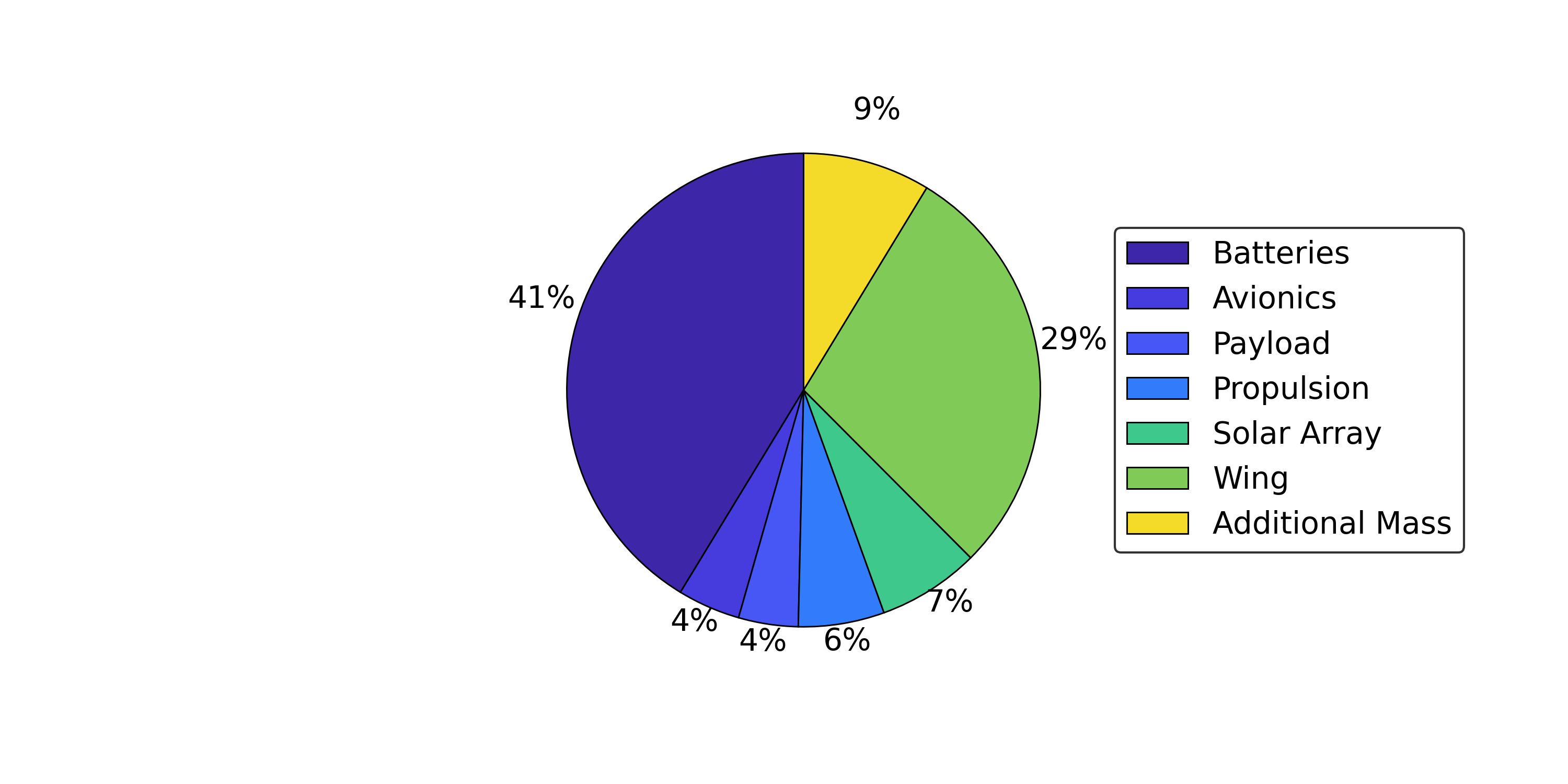}
  \caption{Our case.}
  \label{ecomass}
\end{subfigure}%
\caption{Comparison of the mass breakdowns. The mass is similarly distributed between the different parts of the drone in both cases, even though FBHALE was made using much more complex tools in the design loop.}
\label{compmass}
\end{figure*}

We obtain the convergence graphs of Fig. \ref{convgraph}, and the final design variable values of Tab. \ref{valid}. We add other variables in this table and the values for \cite{colas_hale_2018} as well. Our optimization results are in the same order as those of \cite{colas_hale_2018}. In particular we can see on Fig. \ref{compmass} that the mass is similarly distributed between the different parts of the drone. 

We can see that the share of the weight due to propulsion is nearly the same in both works. This is expected as we derive our propulsion data from \cite{colas_hale_2018} in order to achieve this (see section \ref{propulsion}). However, the shares of the weight due to the batteries, the wing and the solar panels are also similar in both works. This shows that our modelling is a good approximation of this higher-fidelity modelling.
 
\begin{table*}
\centering
\caption{Final design variable values for validation case.}
\label{valid}       % Give a unique label
\begin{tabular}{llll}
\hline\noalign{\smallskip}
Variable & Final values for our case & Final values for FBHALE & Unit\\
\noalign{\smallskip}\hline\noalign{\smallskip}
  Twist control points & [6 12 14 15] & - & deg \\
  Skin thickness control points & [1.8 1.7 2.4 3.5] & - & mm \\
  Spar thickness control points & [2.8 1.1 1.0 1.0] & - & mm \\
  Thickness-to-chord ratio control points & [0.10 0.22 0.23 0.24] & - & - \\
  Span & 53.5 & 45.6 & m \\
  Root chord & 1.4 & - & m \\
  Taper ratio & 0.3 & - & - \\
  Total mass & 245 & 320 & kg \\
  Wing surface & 50.3 & 71.8 & m$^2$\\
  Aspect ratio & 57 & 29 & - \\
  $C_{L}$ at cruise & 1.39 & 1.33 & - \\
  $C_{L}^{3/2}/C_D$ at cruise & 50.2 & 40.1 & - \\
  Motor location over semi-span ratio & 0.33 & 0.46 & - \\
\noalign{\smallskip}\hline
\end{tabular}
\end{table*}

There are however some differences, in particular regarding the aspect ratio. Those are mainly due to our not taking into account 1-cosine gusts in our framework.
Indeed, we only take into account a shear gust wall, whereas \cite{colas_hale_2018} also takes into account a 1-cosine gust, which leads to more stress than a gust wall for a large-span wing.
This difference over-penalize large chords and under-penalize long spans, resulting in a higher aspect-ratio. However, the exact shape of the wing is not the main focus of this work, but rather the coupling between structural mass and total drone mass. This coupling is necessary to see the global influence of the choice of structural materials, and is well captured here. 

A small difference in data can have a big impact on the final mass of the drone. Indeed, a small increase in the drone's total weight, leads to an increase of battery and solar panel weight in order to thrust the heavier drone, and to an increase of wing structural weight to lift the heavier drone. These weight increases contribute to a further increase in the drone's total weight. This is called ``snowball'' effect. For example, the fixed mass of the drone is only 20.5 kg (payload and avionics), but in order to thrust and lift that mass, the final optimized mass of the drone is more than 10 times that.

\subsection{CO2 footprint minimization problem}

Once the framework validated, an optimization is run. It should be recalled that our goal is to optimize the CO2 footprint of the HALE, therefore, the objective function is now this total CO2 emitted. Nevertheless, gradient based solvers converge to a locally optimal point, so the reach of a global optimum depends heavily on the starting point of the optimization. For this reason, a multi-start strategy described in Tab. \ref{multi} is used. Every combination of values is used, meaning that a total of $2 \cdot 3 \cdot 3 \cdot 3 \cdot 4 = 216$ optimizations are run. For each variable, the values are regularly spread between the highest and lowest. For example, for the skin thickness control points, the three values given at the start to that variable are [0.001 0.002 0.003 0.004], [0.0015 0.003 0.0045 0.006], and [0.002 0.004 0.006 0.008].

\begin{table*}
\caption{Design variable starting values for multi-start strategy.}
\label{multi}       % Give a unique label
\begin{tabular}{lllll}
\hline\noalign{\smallskip}
Design variable & Lowest starting value & Highest starting value & \begin{tabular}[c]{@{}l@{}}Number of\\ starting values\end{tabular} & Unit\\
\noalign{\smallskip}\hline\noalign{\smallskip}
  Density & [500 500] & [600 600] & 2 & kg/m$^3$ \\
  Twist control points & [10 15 15 15] & - & 1 & deg \\
  Skin thickness control points & 0.002$\cdot$[0.5 1 1.5 2] & 0.004$\cdot$[0.5 1 1.5 2] & 3 & m \\
  Spar thickness control points & 0.001$\cdot$[1 1 1 1] & 0.003$\cdot$[1 1 1 1] & 3 & m \\
  Thickness-to-chord ratio control points & 0.05$\cdot$[0.75 1 1 1.25] & 0.17$\cdot$[0.75 1 1 1.25] & 3 & - \\
  Span & 25 & 100 & 4 & m \\
  Root chord & 1.5 & - & 1 & m \\
  Taper ratio & 0.3 & - & 1 & - \\
  Motor location over semi-span ratio & 0.3 & - & 1 & - \\
\noalign{\smallskip}\hline
\end{tabular}
\end{table*}

To begin with, the optimization was run with seven possible materials integrated in the material design variable. Their properties appear in Tab. \ref{matos}. Those named material 1, 2 and 3 are roughly homogenized sandwiched panels. Their skins are of the same carbon fibre reinforced polymer (CFRP) as the CFRP material in Tab. \ref{matos}. The expanded polystyrene (PS) foam core of material 1 is replaced by a balsa core in material 2 and a cork core in material 3. Those cores are heavier than the one of material 1 but emit less CO2 per kilogram of material.

Two material indices are computed following the method of \cite{ashby_materials_2004}. Both indices correspond to minimizing the CO2 emitted by the structure $CO2_{struct}$, (defined in Eq. \ref{eqCOstr}), for a flat plate loaded with in-plane compression. The buckling index is useful in the case where the buckling constraint described in section \ref{bucsection} is active. In this case, the index to be maximized is $E^{1/3}/\rho/CO2_{mat}$, where E is the Young’s modulus of the material, $\rho$ its density, and $CO2_{mat}$, the CO2 footprint of the material (defined in Eq. \ref{eqCOmat}). The strength index is useful in the case where the mechanical failure constraint is active. In this case, the index to be maximized is $\sigma_{f}/\rho/CO2_{mat}$, where $\sigma_f$ is the failure strength (yield strength or tensile strength depending on the material). Those material indices appear in Tab. \ref{matos}.
 
\begin{table*}
\caption{Material properties.}
\label{matos}       % Give a unique label
\begin{tabular}{lllllllll}
\hline\noalign{\smallskip}
Property & Material 1 & Material 2 & Material 3 & CFRP & GFRP & Aluminum & Steel & Unit\\
\noalign{\smallskip}\hline\noalign{\smallskip}
  Density & 504.5 & 529 & 560.5 & 1565 & 1860 & 2800 & 7750 & kg/m$^3$ \\
  CO2 emissions & 44.9 & 42.8 & 40.3 & 48.1 & 6.18 & 8.66 & 3.28 & kg$_{CO2}$/kg \\
  Young's modulus & 42.5 & 42.5 & 42.5 & 54.9 & 21.4 & 72.5 & 200 & GPa \\
  Shear modulus & 16.3 &16.3 & 16.3 & 21 & 8.14 & 27 & 78.5 & GPa \\
  Failure strength & 587 & 237 & 587 & 670 & 255 & 445 & 562 & MPa \\
  Buckling index & 0.1539 & 0.1543 & 0.1544 & 0.05049 & 0.24153 & 0.1720 & 0.2302 & N$^{\frac{1}{3}} \cdot$m$^{\frac{7}{3}}$/kg$_{CO2}$ \\
  Strength index & 25885 & 10484 & 25959 & 8901 & 22184 & 18331 & 22119 & N$\cdot$m$\cdot$kg$_{CO2}$ \\
\noalign{\smallskip}\hline
\end{tabular}
\end{table*}

\begin{figure*}
\begin{subfigure}{.49\textwidth}
  \centering
  \includegraphics[width=\linewidth]{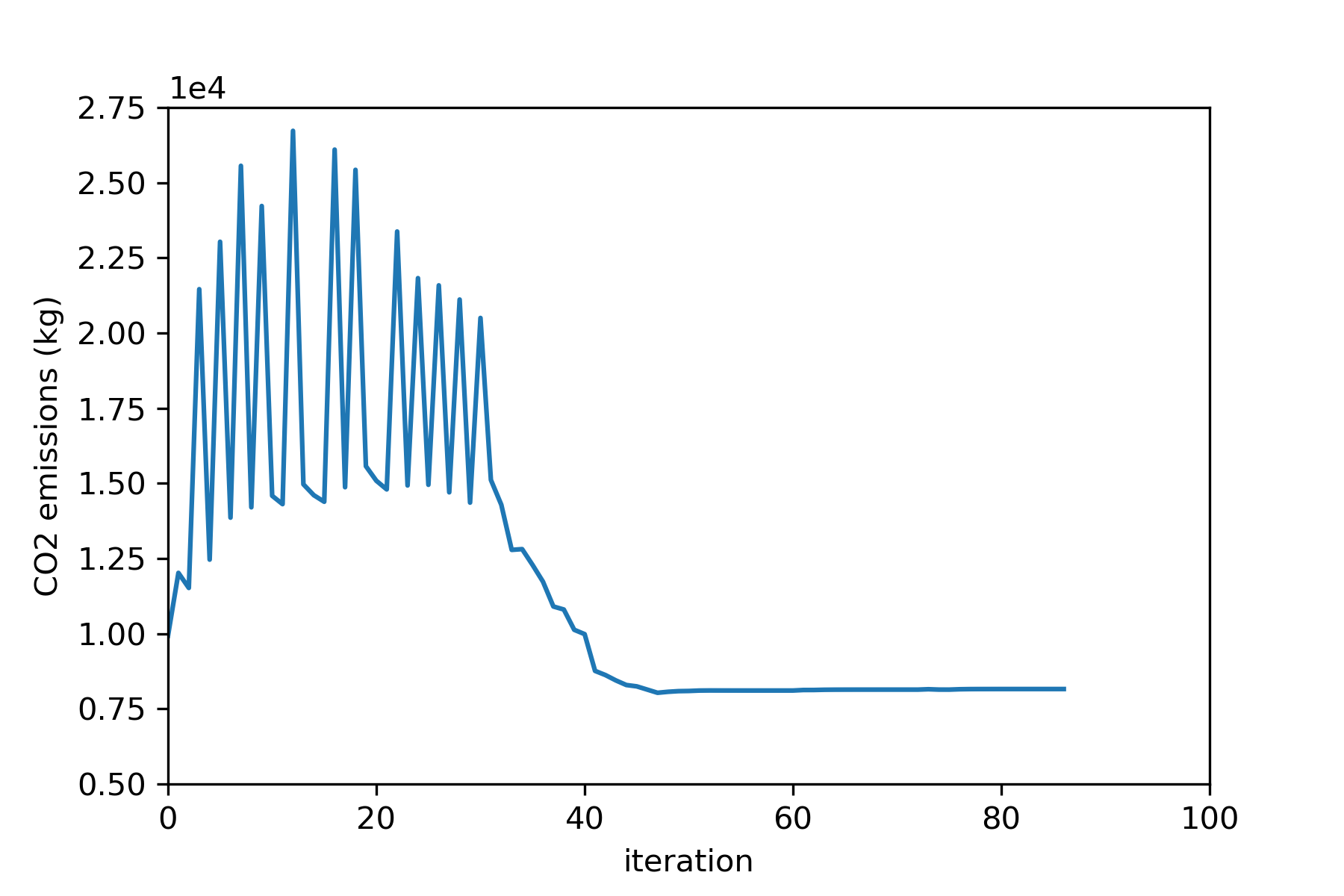}
  \caption{Objective function: total CO2 emitted by the drone.}
  \label{co2}
\end{subfigure}
\begin{subfigure}{.49\textwidth}
  \centering
  \includegraphics[width=\linewidth]{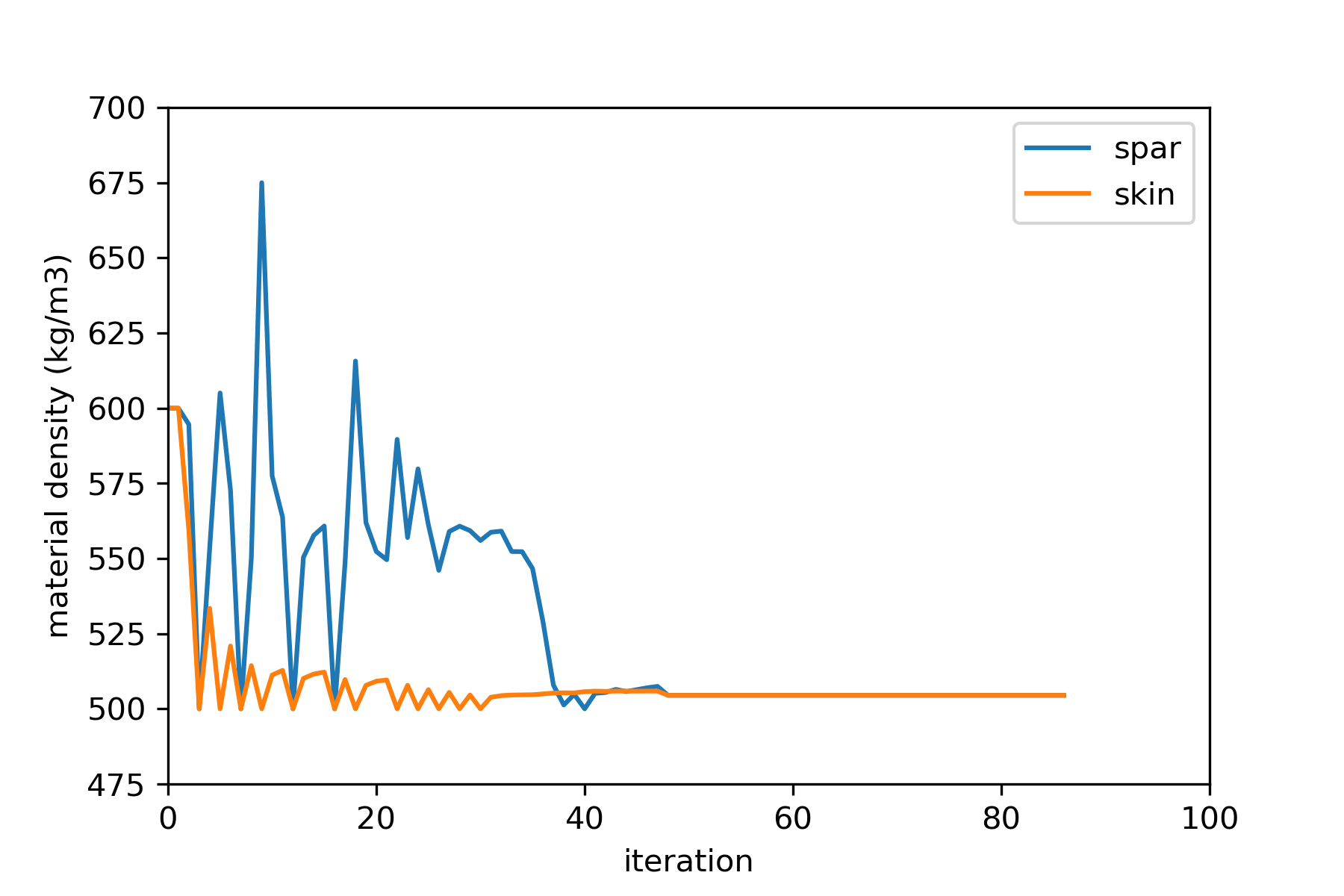}
  \caption{Material density convergence graph.}
  \label{dens}
\end{subfigure}
\caption{Convergence graphs: the lightest material is selected for both spars and skins.}
\label{convgraph2mat}
\end{figure*}

\begin{figure}
  \centering
  \includegraphics[width=0.49\textwidth]{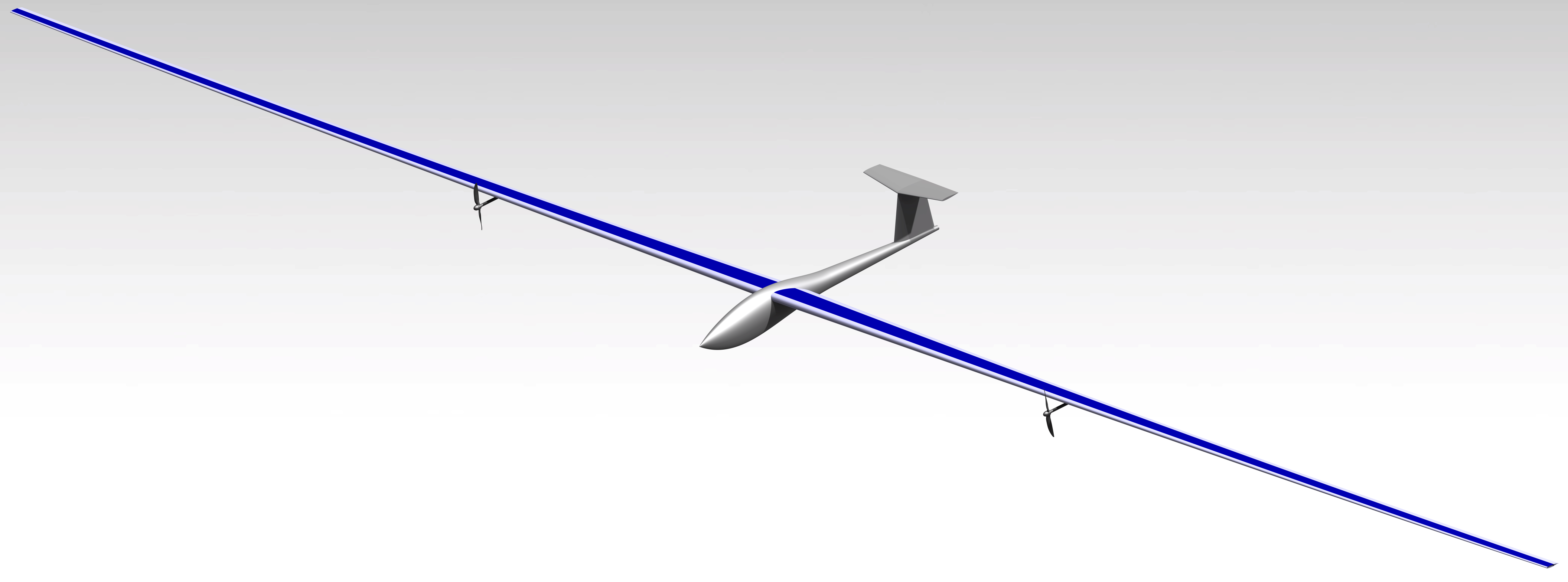}
\caption{Visualization of the CAD model for the optimal HALE wing structure obtained.}
\label{vizu}       % Give a unique label
\end{figure}

The CO2 emissions and density convergence graphs obtained during optimization are shown on Fig. \ref{co2} and \ref{dens}, respectively. In addition, the CAD model of the resulting HALE is shown on Fig. \ref{vizu} (with generic boom and tail shapes since only the wing structure is optimized). The optimization converges towards Material 1 for both spars and skins. This material is the one with the lowest density among those accessible through the material design variable. However, it is not the one with the best material indices. Indeed, the indices for Material 3 are both higher than those of Material 1. A fixed-material optimization is run for both of these materials and confirms that Material 1 preforms better. This material is also the optimal material in case the total weight of the drone is the objective function.

The fact that the CO2 emissions taken into account in the objective function are not only due to the structure (Eq. \ref{eqCOtot}) explains this result. Indeed, although with Material 3 the structure emits less CO2 than with Material 1, it is also heavier. The structure being heavier with Material 3, it needs more power to be thrusted, and thus, more CO2 is emitted by the batteries and solar panels. This compensates the lower emission due to the structure. As a result, material indices as in \cite{ashby_materials_2004} can’t be used in this case. A direct ratio between the weight of the structure and the CO2 emitted by the batteries and solar panels would be needed in order to have a truly useful material index. However, this ratio is not accessible as it evolves during optimization. Therefore, a method as the gradient optimization proposed in this work is necessary in order to choose the best material, if all material candidates can’t be tested individually.

\begin{figure}
  \centering
  \includegraphics[width=0.49\textwidth]{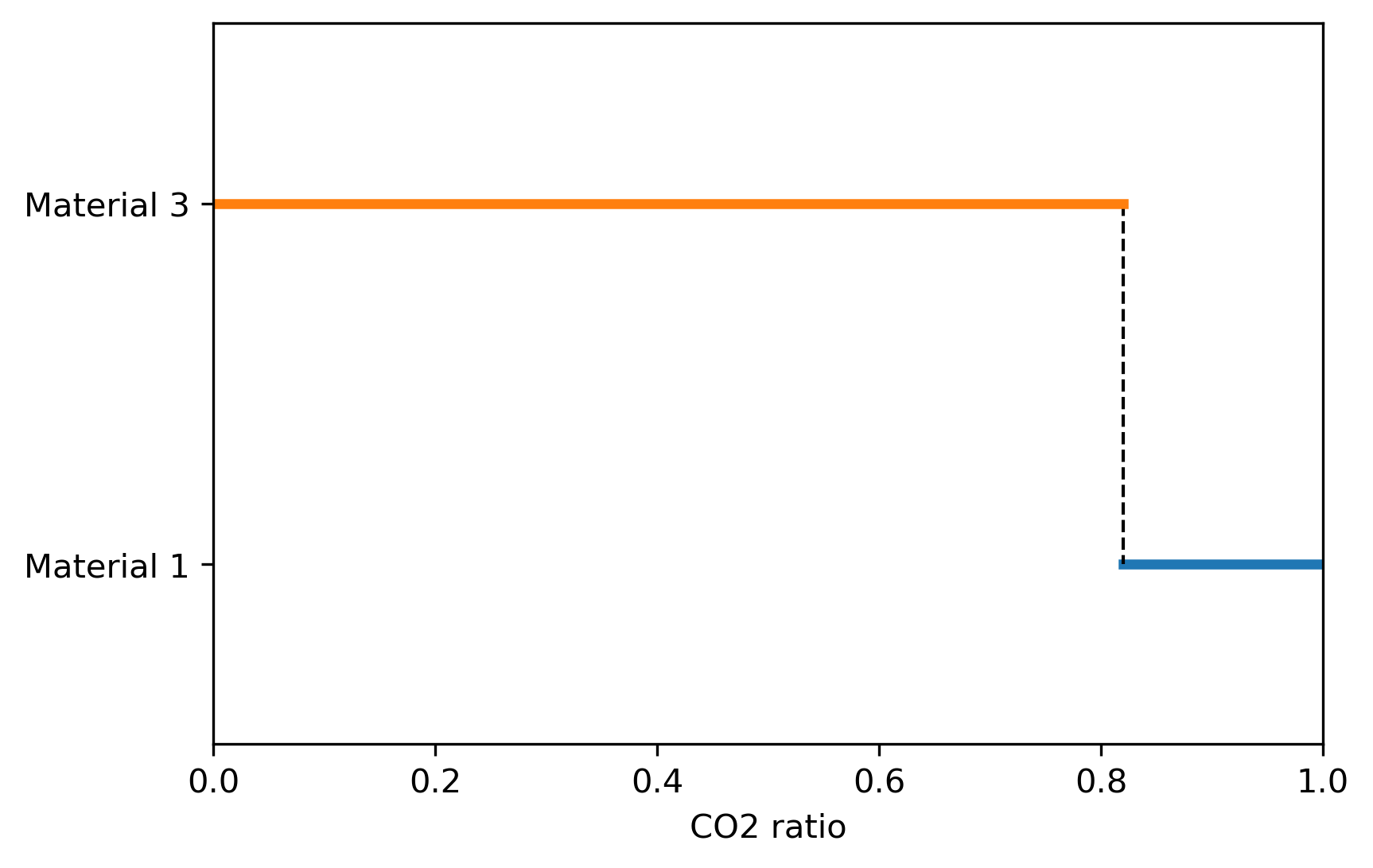}
\caption{Optimal material depending on changes to the CO2 emissions of Material 3. CO2 ratio is the ratio of the CO2 emissions of Material 3 to the CO2 emissions of Material 1.}
\label{optiEmis}       % Give a unique label
\end{figure}

An informal test is carried out to compare the relative influence of the density and the CO2 emissions of a material on its being optimal. For this test, the CO2 emissions of Material 3 are lowered until this material becomes optimal. In this test, Material 1 and Material 3 differ only by their densities and CO2 emissions. The density of Material 3 is 11\% higher than the density of Material 1. The results are presented in Fig. \ref{optiEmis}. It can be seen that the CO2 emissions of Material 3 have to be around 20\% lower than those of material 1, in order to compensate its 11\% higher density. It is reminded that Material 3 is the eco-material surrogate for Material 1, the expanded PS foam core of Material 1 being replaced by a cork core in Material 3. Therefore, in this case, in order to lower the total CO2 emissions of the drone, an eco-material surrogate must also be almost as good as the initial material in terms of weight.

The code (see section \ref{repres}) is not optimized in terms of computational time. This will be addressed in future work. However, the speed of the optimization method presented in this article does not depend on the number of materials in the material catalogue. It is therefore particularly adapted to big catalogues in comparison with a brutal force method consisting in optimizing the structure for every material in the catalogue.

\section{Conclusion}
\label{conclu}

Material choice was successfully integrated in a continuous multidisciplinary design optimization regarding CO2 emissions. The simple HALE model developed showed acceptable agreement with more complex models such as FBHALE, which was designed using much more sophisticated tools.

Our modified version of OpenAeroStruct includes certain improvements with respect to the original version of the tool. For instance, it is adapted to HALE drones and contains more physics: material choice among a discrete catalogue, batteries, solar panels, buckling constraint, etc.

A key result is that a material optimal in terms of drone total weight is also optimal in terms of drone total CO2 emissions, even for an electrical drone. In order to be competitive in terms of total CO2 emissions, an eco-material substitute must be almost as good as the initial in terms of drone total weight. This is due to a ``snowball'' effect on weight.

Finally, a good feature that could be addressed in the future would be a multi-objective optimization between CO2 footprint and cost.

\section*{Replication of results}
\label{repres}
The optimization presented in the results can be obtained by following the instruction of the file howToStart at: \url{https://github.com/mid2SUPAERO/ecoHALE/tree/downloadEcohale}

\section*{Acknowledgements}
The authors would like to thank Mauricio Dwek and Nicolas Martin from Granta Design for their support in using CES Selector.

The authors would also like to thank Fondation ISAE-SUPAERO for funding and enabling the Research Assistant position that lead this article to fruition.


\bibliographystyle{unsrt}  
\bibliography{references}  %%% Remove comment to use the external .bib file (using bibtex).

%

\end{document}